# MODERATE DEVIATIONS FOR POISSON–DIRICHLET DISTRIBUTION

BY SHUI FENG[1] AND FUQING GAO[2]

*McMaster University and Wuhan University*

The Poisson–Dirichlet distribution arises in many different areas. The parameter $\theta$ in the distribution is the scaled mutation rate of a population in the context of population genetics. The limiting case of $\theta$ approaching infinity is practically motivated and has led to new, interesting mathematical structures. Laws of large numbers, fluctuation theorems and large-deviation results have been established. In this paper, moderate-deviation principles are established for the Poisson–Dirichlet distribution, the GEM distribution, the homozygosity, and the Dirichlet process when the parameter $\theta$ approaches infinity. These results, combined with earlier work, not only provide a relatively complete picture of the asymptotic behavior of the Poisson–Dirichlet distribution for large $\theta$, but also lead to a better understanding of the large deviation problem associated with the scaled homozygosity. They also reveal some new structures that are not observed in existing large-deviation results.

**1. Introduction.** For $\theta > 0$, let $\sigma_1(\theta) \geq \sigma_2(\theta) \geq \cdots$ be the points of a nonhomogeneous Poisson process with mean measure density

$$\theta u^{-1} e^{-u}, \qquad u > 0,$$

and $\sigma(\theta) = \sum_{i=1}^{\infty} \sigma_i(\theta)$. Set

$$(1.1) \qquad \mathbf{P}(\theta) = (P_1(\theta), P_2(\theta), \ldots) = \left(\frac{\sigma_1(\theta)}{\sigma(\theta)}, \frac{\sigma_2(\theta)}{\sigma(\theta)}, \ldots\right).$$

Then it is known that $\mathbf{P}(\theta)$ and $\sigma(\theta)$ are independent, and $\sigma(\theta)$ is a Gamma$(\theta, 1)$-distributed random variable. The law of $\mathbf{P}(\theta)$ is called the Poisson–Dirichlet distribution with parameter $\theta$, and is denoted by PD$(\theta)$.

Received May 2007; revised October 2007.
[1]Supported by the Natural Science and Engineering Research Council of Canada.
[2]Supported by the NSF of China (No. 10571139).

Let $U_k, k = 1, 2, \ldots$, be a sequence of i.i.d. random variables with common distribution, Beta$(1, \theta)$. Set

$$(1.2) \quad X_1(\theta) = U_1, \qquad X_n(\theta) = (1 - U_1) \cdots (1 - U_{n-1}) U_n, \qquad n \geq 2.$$

Then with probability one

$$\sum_{k=1}^{\infty} X_k(\theta) = 1,$$

and the law of $(X_1(\theta), X_2(\theta), \ldots)$ is called the GEM distribution, denoted by GEM$(\theta)$. The law of the descending order statistics $X_{(1)}(\theta) \geq X_{(2)}(\theta) \geq \cdots$ of $X_1(\theta), X_2(\theta), \ldots$ is also PD$(\theta)$.

Let $\xi_k, k = 1, \ldots$, be a sequence of i.i.d. random variables, independent of $\mathbf{P}(\theta)$, with a common diffusive distribution $\nu$ on $[0, 1]$, that is, $\nu(\{x\}) = 0$ for every $x$ in $[0, 1]$. Set

$$(1.3) \qquad \Xi_{\theta,\nu} = \sum_{k=1}^{\infty} P_k(\theta) \delta_{\xi_k}.$$

We call the law of $\Xi_{\theta,\nu}$, the Dirichlet process, denoted by Dirichlet$(\theta, \nu)$.

The Poisson–Dirichlet distribution was introduced by Kingman [21] to describe the distribution of gene frequencies in a large neutral population at a particular locus. The component $P_k(\theta)$ represents the proportion of the $k$th most frequent allele. If $\varepsilon$ is the individual mutation rate and $N_e$ is the effective population size, then the parameter $\theta = 4N_e\varepsilon$ is the scaled population mutation rate. The GEM distribution can be obtained from the Poisson–Dirichlet distribution through a procedure called *size-biased sampling*. It provides an effective way of doing calculations involving the Poisson–Dirichlet distribution. The name, GEM distribution, was coined by Ewens after Grifffiths, Engen and McCloskey for their contributions to the development of the structure. The Dirichlet process first appeared in [11] in the context of Bayesian statistics. It can be viewed as a labelled version of the Poisson–Dirichlet distribution. More background information can be found in [8].

For any integer $m \geq 2$, consider a random sample of size $m$ from a population following the Poisson–Dirichlet distribution. Given the population proportion, $\mathbf{p} = (p_1, p_2, \ldots)$, the probability that all samples are of the same type is given by

$$H_m(\mathbf{p}) = \sum_{i=1}^{\infty} p_i^m.$$

The quantity $H_2(\mathbf{p})$ is called the population homozygosity. It is an important statistic in population genetics. For general $m$, we refer to $H_m(\mathbf{p})$, as the homozygosity of order $m$.

Consider a family of random variables $\{Y_\lambda : \lambda > 0\}$. Assume a law of large numbers holds; that is, $Y_\lambda$ converges in distribution to a constant $c$ as $\lambda$ approaches infinity. A fluctuation theorem such as the central limit theorem is a statement that there exists a function $b(\lambda)$ approaching infinity for large $\lambda$ such that

$$b(\lambda)(Y_\lambda - c) \Rightarrow Y, \qquad \lambda \to \infty,$$

where $Y$ is a nontrivial random variable and "$\Rightarrow$" denotes convergence in distribution. A large-deviation result is concerned with estimates of probabilities $P\{Y_\lambda - c \in A\}$ for measurable sets $A$. A moderate-deviation result lies between the fluctuation theorems and large deviations. It is concerned with estimates of probabilities $P\{a(\lambda)(Y_\lambda - c) \in A\}$ for measurable sets $A$, where $a(\lambda)$ is an intermediate scale between 1 and $b(\lambda)$.

The objective of this paper is to establish moderate-deviation principles (henceforth MDP) for $\text{GEM}(\theta)$, $\text{PD}(\theta)$, the homozygosity and $\text{Dirichlet}(\theta, \nu)$, when $\theta$ approaches infinity.

The study of the behavior of $\mathbf{P}(\theta) = (P_1(\theta), P_2(\theta), \ldots)$ for large $\theta$, goes back to the seventies. In Watterson and Guess [29], $E[P_1(\theta)]$ was shown to be asymptotically $\log\theta/\theta$. Griffiths [15] obtained the explicit weak limit of $\theta\mathbf{P}(\theta)$ and a central limit theorem for the population homozygosity. The limiting case of large $\theta$ is equivalent to a situation where the mutation rate per individual is fixed and the effective population size is large. Motivated by the work of Gillespie [12] on the role of population size in molecular evolution, there have been renewed interests in the asymptotic behavior of $\text{PD}(\theta)$ for large $\theta$ (see [4, 9, 18, 19, 20]). In particular, in [19], central limit theorems are obtained for the homozygosity of order $m$. Large deviations are established in [4] for $\text{PD}(\theta)$ and the homozygosity, and in [9] for the GEM distribution. Large deviations for $\text{Dirichlet}(\theta, \nu)$ can be found in [3, 23].

Although MDP is a natural mathematical object that warrants rigorous investigation and our study does reveal some new mathematical phenomena, the real motivation for this work comes from the results in [4, 19]. It was shown in [19] that, as $\theta$ goes to infinity,

$$(1.4) \qquad \frac{\theta^{m-1}}{\Gamma(m)} H_m(\mathbf{P}(\theta)) \to 1$$

and

$$(1.5) \qquad \sqrt{\theta}\left(\frac{\theta^{m-1}}{\Gamma(m)} H_m(\mathbf{P}(\theta)) - 1\right) \to \mathcal{Z}(m),$$

where $\mathcal{Z}(m)$ is a normal random variable with mean zero and variance $\frac{\Gamma(2m)}{\Gamma^2(m)} - m^2$. These are the law of large numbers and central limit theorem for $\frac{\theta^{m-1}}{\Gamma(m)} H_m(\mathbf{P}(\theta))$. A natural companion to these limit theorems is the large deviations $\frac{\theta^{m-1}}{\Gamma(m)}H_m(\mathbf{P}(\theta))$ from one, or equivalently the large deviations of $\frac{\theta^{m-1}}{\Gamma(m)}H_m(\mathbf{P}(\theta)) - 1$ from zero. Unfortunately this problem is still open. The large deviation principle established in [4] is concerned with the deviations of $H_m(\mathbf{P}(\theta))$ from zero. The scale difference between $H_m(\mathbf{P}(\theta))$ and $\frac{\theta^{m-1}}{\Gamma(m)}H_m(\mathbf{P}(\theta))$ is of order of $\theta^{m-1}$. Multiplying $\frac{\theta^{m-1}}{\Gamma(m)}H_m(\mathbf{P}(\theta)) - 1$ by a factor $\theta^\gamma$, places us in the territory of MDP. One would hope that the study of MDP will shed light on resolving the large deviation problem which corresponds to $\gamma = 0$. The MDPs we obtain require that $\gamma$ is bigger than a strictly positive number. Thus a gap exists between the MDPs and the LDP. This seems to indicate that a large deviation principle may not exist for $\frac{\theta^{m-1}}{\Gamma(m)}H_m(\mathbf{P}(\theta)) - 1$.

This paper is organized as follows. The basic terminology of LDP, MDP and a comparison lemma are given in Section 2. In Section 3, we discuss the fluctuation theorems associated with $PD(\theta)$, Dirichlet$(\theta, \nu)$ and the homozygosity. A new proof is given for the central limit theorem of homozygosity in [19], using Campbell's theorem. A MDP for GEM is established in Section 4. Section 5 deals with the MDP for $PD(\theta)$. Since the condition of the Gärtner–Ellis theorem is not satisfied, we prove the result by direct calculation. The MDP obtained in Section 6, is for the homozygosity, for which the MDP holds in a narrower range of scales. The proof is based mainly on Campbell's theorem. In the MDP literature, general results such as those in [7, 13, 30], usually require the finiteness of exponential moments in a small neighborhood of zero so that the Laplace method can be used. Here the exponential moment is infinite on the positive half-line. One way to deal with the infinite exponential moment is to verify Ledoux's condition in [22]. Since this does not seem easy to do, we choose the truncation method instead. Finally in Section 7, we establish the MDP for Dirichlet$(\theta, \nu)$. Compared to the Sanov theorem, the LDP rate function for Dirichlet$(\theta, \nu)$ is a reversed form of relative entropy. Here the MDP rate function for Dirichlet process is the same as the MDP rate function for the empirical process of an i.i.d. random sequence with common distribution $\nu$. When $\nu$ is supported on a finite number of points, one can see this clearly from the fact that both the relative entropy and its reversed form have the same second-order derivative at $\nu$.

The MDPs for the Poisson–Dirichlet distribution and GEM have a different speed from the MDPs for the homozygosity and the Dirichlet process, the latter having a more standard structure. One explanation for this is that in the cases of the Poisson–Dirichlet distribution and GEM, we are concerned with partial information such as alleles with a certain proportion size or age order, while for the homozygosity and the Dirichlet process, all alleles contribute. One expects that similar results and structures exist for the two-parameter Poisson–Dirichlet distribution and Dirichlet process [10, 26].

**2. Preliminaries.** In this section we introduce the terminology on LDP and MDP used in this paper, and prove a comparison lemma that plays an important role in proving the main results. Comprehensive coverage on LDP techniques can be found in [6].

DEFINITION 2.1. Let $E$ be a Polish space with metric $d$, and $\{Y_\theta : \theta > 0\}$ be a family of $E$-valued random variables. Denote the law of $Y_\theta$ by $P_\theta$.

(1) The family of probability measures $\{P_\theta : \theta > 0\}$ (or the family $\{Y_\theta : \theta > 0\}$) is said to satisfy a LDP with speed $\lambda(\theta)$ and rate function $I(\cdot)$, if for any closed set $F$ and open set $G$ in $E$

$$\limsup_{\theta \to \infty} \lambda(\theta) \log P_\theta\{F\} \leq -\inf_{x \in F} I(x),$$

$$\liminf_{\theta \to \infty} \lambda(\theta) \log P_\theta\{G\} \geq -\inf_{x \in G} I(x),$$

for any $c > 0$, $\{x : I(x) \leq c\}$ is compact.

In short form, we say $(P_\theta, I(\cdot), \lambda(\theta))$ satisfies a LDP.

(2) The family $\{P_\theta : \theta > 0\}$ is said to satisfy a local LDP with speed $\lambda(\theta)$ and rate function $I(\cdot)$, if for every $x$ in $E$

$$\lim_{\delta \to 0} \limsup_{\theta \to \infty} \lambda(\theta) \log P\{d(Y_\theta, x) \leq \delta\}$$
$$= \lim_{\delta \to 0} \liminf_{\theta \to \infty} \lambda(\theta) \log P\{d(Y_\theta, x) < \delta\} = -I(x),$$

and for any $c > 0$, $\{x : I(x) \leq c\}$ is compact.

(3) The family $\{P_\theta : \theta > 0\}$ is exponentially tight with speed $\lambda(\theta)$ if for every $L > 0$, there is a compact set $K_L$ in $E$ such that

$$\limsup_{\theta \to \infty} \lambda(\theta) \log P\{Y_\theta \notin K_L\} \leq -L.$$

REMARK 2.1. It is known that a local LDP combined with exponential tightness implies the LDP (cf. [27]).

DEFINITION 2.2. We use $\Rightarrow$ to denote convergence in distribution.

(1) The family $\{Y_\theta : \theta > 0\}$ is said to satisfy a fluctuation theorem if there exist functions $b(\theta)$, $c(\theta)$ and a finite nondeterministic random variable $Z$ such that

$$\lim_{\theta \to \infty} b(\theta) = \infty, \qquad b(\theta)[Y_\theta - c(\theta)] \Rightarrow Z, \qquad \theta \to \infty.$$

(2) Assume that the family $\{Y_\theta : \theta > 0\}$ satisfies the fluctuation theorem above. Let $a(\theta)$ satisfy

(2.1) $$\lim_{\theta \to \infty} a(\theta) = \infty, \qquad \lim_{\theta \to \infty} \frac{a(\theta)}{b(\theta)} = 0.$$

The family $\{P_\theta : \theta > 0\}$ or equivalently the family $\{Y_\theta : \theta > 0\}$ is said to satisfy a MDP with speed $\lambda(\theta)$ [depending on $a(\theta)$] and rate function $I(\cdot)$ if the family $\{a(\theta)[Y_\theta - c(\theta)] : \theta > 0\}$ satisfies a LDP with speed $\lambda(\theta)$ and rate function $I(\cdot)$. Thus the MDP for $\{Y_\theta : \theta > 0\}$ is the LDP for $\{a(\theta)[Y_\theta - c(\theta)] : \theta > 0\}$.

The next lemma is a useful tool in deriving the MDPs of this paper.

LEMMA 2.1. *Let $\{\xi_\theta : \theta > 0\}$ and $\{\eta_\theta > 0 : \theta > 0\}$ be two families of real-valued random variables. Assume that for any $\delta > 0$*

(2.2) $$\limsup_{\theta \to \infty} \lambda(\theta) \log P(|\eta_\theta - 1| \geq \delta) = -\infty.$$

*Then $(P(\xi_\theta \in \cdot), I(x), \lambda(\theta))$ satisfies a LDP iff $(P(\xi_\theta \eta_\theta \in \cdot), I(x), \lambda(\theta))$ satisfies a LDP.*

PROOF. For any $\delta > 0$, choose $\tilde{\delta} = \min\{\frac{\delta}{2}, \frac{1}{2}\}$. Then it is clear that

$$\{|\eta_\theta^{-1} - 1| \geq \delta\} \subset \{|\eta_\theta - 1| \geq \tilde{\delta}\},$$

which, combined with (2.2), implies

(2.3) $$\limsup_{\theta \to \infty} \lambda(\theta) \log P(|\eta_\theta^{-1} - 1| \geq \delta) = -\infty.$$

For any $x \in \mathbb{R}$, $\gamma > 0$ and $\delta > 0$,

$$P(|\xi_\theta \eta_\theta - x| \leq \gamma)$$
$$\leq P(|\eta_\theta^{-1} - 1| \geq \delta) + P(|\eta_\theta^{-1} - 1| \leq \delta, |\xi_\theta - x| \leq \gamma \eta_\theta^{-1} + |x||\eta_\theta^{-1} - 1|)$$
$$\leq P(|\eta_\theta^{-1} - 1| \geq \delta) + P(|\xi_\theta - x| \leq \gamma(1 + \delta) + |x|\delta)$$

which implies

$$\lim \binom{\sup}{\inf} \lambda(\theta) \log P(|\xi_\theta \eta_\theta - x| \leq \gamma) \leq \lim \binom{\sup}{\inf} \lambda(\theta) \log P(|\xi_\theta - x| \leq \gamma).$$

Symmetrically,

$$\lim \binom{\sup}{\inf} \lambda(\theta) \log P(|\xi_\theta - x| \leq \gamma) \leq \lim \binom{\sup}{\inf} \lambda(\theta) \log P(|\xi_\theta \eta_\theta - x| \leq \gamma).$$

Furthermore, for any $L > 0$,

$$P(|\xi_\theta \eta_\theta| \geq L) \leq P(|\eta_\theta - 1| \geq \delta) + P(|\xi_\theta| \geq (1+\delta)^{-1} L),$$
$$P(|\xi_\theta| \geq L) \leq P(|\eta_\theta^{-1} - 1| \geq \delta) + P(|\xi_\theta \eta_\theta| \geq (1+\delta)^{-1} L).$$

Thus the exponential tightness of $\{\xi_\theta : \theta > 0\}$ is equivalent to the exponential tightness of $\{\xi_\theta \eta_\theta : \theta > 0\}$. The lemma now follows from Remark 2.1.
□

**3. Fluctuation theorems.** We start this section with a discussion of the asymptotic behavior of the random variable $\sigma(\theta)$ for large $\theta$. It plays a key role in connecting the Poisson process to the Poisson–Dirichlet distribution. To put our MDP results into perspective, we present in this section several known fluctuation theorems for the Poisson–Dirichlet distribution, the Dirichlet process and the homozygosity of order $m$.

3.1. *Gamma distribution.* Recall that $\sigma(\theta)$ is a Gamma$(\theta, 1)$ random variable with density

$$(3.1) \qquad \frac{1}{\Gamma(\theta)} u^{\theta-1} e^{-u}, \qquad 0 < u < \infty,$$

and exponential moment

$$(3.2) \qquad E[e^{t\sigma(\theta)}] = \begin{cases} \frac{1}{(1-t)^\theta}, & t < 1 \\ \infty, & \text{else.} \end{cases}$$

Let

$$(3.3) \qquad \Lambda(t) = \begin{cases} \log \frac{1}{1-t}, & t < 1 \\ \infty, & \text{else.} \end{cases}$$

Routine calculations and Gärtner–Ellis theorem lead to the following theorem:

THEOREM 3.1. *When $\theta$ approaches infinity, the following hold:*

(a) $\lim_{\theta \to \infty} \frac{\sigma(\theta)}{\theta} = 1$.
(b) $\sqrt{\theta}(\frac{\sigma(\theta)}{\theta} - 1) \Rightarrow Z$, *where $Z$ is a standard normal random variable.*
(c) *The family of the laws of $\sigma(\theta)/\theta$ satisfies a LDP with speed $1/\theta$ and rate function*

$$(3.4) \qquad I(u) = \begin{cases} u - 1 - \log u, & u > 0 \\ \infty, & \text{else.} \end{cases}$$

Let $a(\theta)$ be a positive function satisfying

(3.5) $$\lim_{\theta \to \infty} \frac{a(\theta)}{\theta} = 0, \qquad \lim_{\theta \to \infty} a(\theta) = \infty.$$

COROLLARY 3.1. *For each $\delta > 0$,*

(3.6) $$\limsup_{\theta \to \infty} \frac{a(\theta)}{\theta} \log P\left\{ \left| \frac{\sigma(\theta)}{\theta} - 1 \right| > \delta \right\} = -\infty$$

*and*

(3.7) $$\limsup_{\theta \to \infty} \frac{a(\theta)}{\theta} \log P\left\{ \left| \frac{\theta}{\sigma(\theta)} - 1 \right| > \delta \right\} = -\infty.$$

PROOF. Equality (3.6) is derived directly from Theorem 3.1. Since

$$P\left\{ \left| \frac{\theta}{\sigma(\theta)} - 1 \right| > \delta \right\} \leq P\left\{ \left| \frac{\sigma(\theta)}{\theta} - 1 \right| > \delta/3 \right\} + P\left\{ \left| \frac{\sigma(\theta)}{\theta} - 1 \right| > 1/2 \right\},$$

one gets (3.7) from (3.6). □

Let $a(\theta)$ be a positive function satisfying

(3.8) $$\lim_{\theta \to \infty} \frac{a(\theta)}{\sqrt{\theta}} = 0, \qquad \lim_{\theta \to \infty} a(\theta) = \infty.$$

The following theorem is standard.

THEOREM 3.2. *The family of the laws of $\sigma(\theta)/\theta$ satisfies a MDP with speed $a^2(\theta)/\theta$ and rate function*

(3.9) $$S(u) = \frac{u^2}{2}, \qquad -\infty < u < \infty.$$

3.2. *Fluctuations.* Consider a nonhomogeneous Poisson process with mean measure

$$e^{-u} du, \qquad -\infty < u < +\infty.$$

Let $\zeta_1 \geq \zeta_2 \geq \cdots$ be the sequence of the points of the nonhomogeneous Poisson process in descending order. Then for each $r \geq 1$ the joint density of $(\zeta_1, \ldots, \zeta_r)$ is

(3.10) $$e^{-\sum_{i=1}^r u_i} e^{-e^{-u_r}}, \qquad -\infty < u_r < \cdots < u_1 < \infty.$$

Let $\beta(\theta) = \log \theta - \log \log \theta$. The following result is obtained in [15].

THEOREM 3.3. *The sequence* $(\theta P_1(\theta) - \beta(\theta), \theta P_2(\theta) - \beta(\theta), \ldots)$ *converges to* $(\zeta_1, \zeta_2, \ldots)$ *in distribution as* $\theta$ *tends to infinity.*

The next theorem is obtained in [19]. We give a different proof here using Campbell's theorem.

THEOREM 3.4. *Let*

$$A_k(\theta) = \sqrt{\theta}\left(\frac{\theta^{k-1}}{\Gamma(k)} H_k(\mathbf{P}(\theta)) - 1\right), \qquad k = 2, 3, \ldots$$

*and* $\mathbf{A}_\theta = (A_2(\theta), A_3(\theta), \ldots)$. *Then*

(3.11) $$\mathbf{A}_\theta \Rightarrow \mathbf{A}, \qquad \theta \to \infty,$$

*where* $\mathbf{A}$ *is a* $\mathbb{R}^\infty$-*valued random element and for each* $r \geq 2$, $(A_2, \ldots, A_r)$ *has a multivariate normal distribution with zero mean and covariance matrix*

(3.12) $$\operatorname{Cov}(A_k, A_l) = \frac{\Gamma(k+l) - \Gamma(k+1)\Gamma(l+1)}{\Gamma(k)\Gamma(l)}, \qquad k, l = 2, \ldots, r.$$

PROOF. For each $k \geq 1$, set

$$B_k(\theta) = \sqrt{\theta}\left(\frac{1}{\Gamma(k)\theta} \sum_{l=1}^\infty \sigma_l^k(\theta) - 1\right),$$

$$\mathbf{B}_\theta = (B_1(\theta), \ldots).$$

For each fixed $r \geq 1$ and any $(\alpha_1, \ldots, \alpha_r)$ in $R^r$, set

$$f(x) = \sum_{k=1}^r \frac{1}{\Gamma(k)\sqrt{\theta}} \alpha_k x^k.$$

It follows from Campbell's theorem that

$$E(e^{[it\sum_{k=1}^r \alpha_k B_k(\theta)]})$$
$$= e^{[-it\sum_{k=1}^r \alpha_k \sqrt{\theta}]} E(e^{[it\sum_{l=1}^\infty f(\sigma_l(\theta))]})$$
(3.13) $$= e^{[-it\sum_{k=1}^r \alpha_k \sqrt{\theta}]} \exp\left\{\theta \int_0^\infty (e^{itf(y)} - 1) y^{-1} e^{-y}\, dy\right\}$$
$$\to \exp\left\{-\frac{t^2}{2} \sum_{j,k=1}^r \alpha_j \alpha_k \frac{\Gamma(j+k)}{\Gamma(j)\Gamma(k)}\right\}.$$

Let $\mathbf{B} = (B_1, \ldots)$ be such that for each $r \geq 1$, $(B_1, \ldots, B_r)$ is a multivariate normal random vector with zero mean and covariance matrix

(3.14) $$\frac{\Gamma(j+k)}{\Gamma(j)\Gamma(k)}, \qquad j, k = 1, \ldots, r.$$

Then (3.13) implies that $\mathbf{B}_\theta$ converges in distribution to $\mathbf{B}$.

For $r \geq 2$, it follows from (1.1) that the following relation holds between $(A_2(\theta), \ldots, A_r(\theta))$ and $(B_2(\theta), \ldots, B_r(\theta))$:

$$(3.15) \qquad A_k(\theta) = B_k(\theta) + \sqrt{\theta}\left(\left(\frac{\theta}{\sigma(\theta)}\right)^k - 1\right)\left(\sum_{l=1}^\infty \frac{\sigma_l^k(\theta)}{\Gamma(k)\theta}\right).$$

It follows from the convergence of $\mathbf{B}_\theta$ to $\mathbf{B}$ that

$$(3.16) \qquad \sum_{l=1}^\infty \frac{\sigma_l^k(\theta)}{\Gamma(k)\theta} \to 1 \quad \text{in distribution.}$$

By Theorem 3.1 and basic algebra, one gets

$$(3.17) \qquad \sqrt{\theta}\left(\left(\frac{\theta}{\sigma(\theta)}\right)^k - 1\right) \Rightarrow -kB_1.$$

By (3.15)–(3.17), one gets

$$(3.18) \qquad \sum_{k=2}^r \alpha_k A_k(\theta) \Rightarrow \sum_{k=2}^r \alpha_k(B_k - kB_1).$$

The theorem now follows from the fact that the covariance of $(B_k - kB_1)$ and $(B_l - lB_1)$ is $\frac{\Gamma(k+l)-\Gamma(k+1)\Gamma(l+1)}{\Gamma(k)\Gamma(l)}$. □

Let $\{X(t), t \in [0, \infty)\}$ be a Gamma process; that is, a stochastic process with stationary independent increments and right-continuous paths with $X(0) = 0$ and such that $X(1)$ has an exponential distribution with parameter 1. For each Borel measurable set $A$, define

$$X_{\theta,\nu}(A) = X(\theta\nu(A))$$

and

$$Z_{\theta,\nu}(A) = \frac{X(\theta\nu(A))}{X(\theta)},$$

where $\nu$ is a diffusive distribution on $[0,1]$, that is, $\nu(\{x\}) = 0$ for every $x$ in $[0,1]$. Set $\nu(t) = \nu([0,t])$, $X_{\theta,\nu}(t) = X_{\theta,\nu}([0,t])$ and $Z_{\theta,\nu}(t) = Z_{\theta,\nu}([0,t])$. Then $Z_{\theta,\nu}(\cdot)$, as a random measure, is distributed as Dirichlet$(\theta, \nu)$. Let $D([0,1])$ be the space of all real-valued cadlag functions defined on $[0,1]$ that are left continuous at 1, equipped with the topology of uniform convergence. Then the functional central limit theorem for processes with independent increments yields immediately that $(X_{\theta,\nu}(t) - \theta\nu(t))/\sqrt{\theta}$ converges to $B(\nu(t))$ in distribution on $D([0,1])$, where $B(t)$ is a standard Brownian motion (cf. VII 3.5 in [16], page 373). This, combined with the fact that $X(\theta)/\theta$ converges to 1, implies the following result:

THEOREM 3.5. *The family of processes $\{\sqrt{\theta}(Z_{\theta,\nu}(t) - \nu(t)), \theta > 0\}$ converges to $B(\nu(t))$ in distribution on $D([0,1])$.*

**4. Moderate deviations for GEM.** Let $a(\theta)$ satisfy (3.5).

The MDP for GEM is thus the LDP for the family of $\{a(\theta)(X_1(\theta), X_2(\theta), \ldots):\theta > 0\}$ when $\theta$ approaches infinity. The result is proved through explicit calculations.

THEOREM 4.1. *The family* $\{P(a(\theta)(X_1(\theta), X_2(\theta), \ldots) \in \cdot) : \theta > 0\}$ *satisfies a LDP on $\mathbb{R}_+^\infty$ with speed $\frac{a(\theta)}{\theta}$ and rate function*

$$I(x_1, x_2, \ldots) = \sum_{i=1}^\infty x_i. \tag{4.1}$$

PROOF. Let us first prove the local LDP. For any $\mathbf{x}, \mathbf{y}$ in $\mathbb{R}_+^\infty$, set

$$|\mathbf{x} - \mathbf{y}| = \sum_{i=1}^\infty \frac{|x_i - y_i| \wedge 1}{2^i}.$$

For any $\mathbf{x}$ in $\mathbb{R}_+^\infty$ and any $\delta > 0$, one can choose $n$ sufficiently large that $\sum_{i=n}^\infty \frac{1}{2^i} < \delta/2$. Then for $\delta_1 < \delta/2$, we have

$$\left\{\mathbf{y} \in \mathbb{R}_+^\infty : \max_{1 \leq i \leq n} |y_i - x_i| < \delta_1\right\} \subset \{\mathbf{y} \in \mathbb{R}_+^\infty : |\mathbf{y} - \mathbf{x}| < \delta\}. \tag{4.2}$$

By taking limits in (4.2), in the order $\theta \to \infty, \delta_1 \to 0, n \to \infty, \delta \to 0$, it follows that

$$\begin{aligned}
&\lim_{\delta \to 0} \liminf_{\theta \to \infty} \frac{a(\theta)}{\theta} \log P(|a(\theta)(X_1(\theta), X_2(\theta), \ldots) - \mathbf{x}| < \delta) \\
&\geq \lim_{n \to \infty} \lim_{\delta_1 \to 0} \liminf_{\theta \to \infty} \frac{a(\theta)}{\theta} \log P\left(\max_{1 \leq i \leq n} |a(\theta)X_i(\theta) - x_i| < \delta_1\right).
\end{aligned} \tag{4.3}$$

On the other hand, for any $m \geq 1$ such that $2^{-m} > \delta$

$$\{\mathbf{y} \in \mathbb{R}_+^\infty : |\mathbf{y} - \mathbf{x}| \leq \delta\} \subset \left\{\mathbf{y} \in \mathbb{R}_+^\infty : \max_{1 \leq i \leq m} |y_i - x_i| \leq 2^m \delta\right\}.$$

Thus for any $\delta_2 < \delta$,

$$\{\mathbf{y} \in \mathbb{R}_+^\infty : |\mathbf{y} - \mathbf{x}| \leq \delta_2\} \subset \left\{\mathbf{y} \in \mathbb{R}_+^\infty : \max_{1 \leq i \leq m} |y_i - x_i| \leq 2^m \delta\right\}. \tag{4.4}$$

By taking the limits in (4.4), in the order $\theta \to \infty$, $\delta_2 \to 0$, $\delta \to 0$, $m \to \infty$, it follows that

$$\begin{aligned}
&\lim_{\delta_2 \to 0} \limsup_{\theta \to \infty} \frac{a(\theta)}{\theta} \log P(|a(\theta)(X_1(\theta), X_2(\theta), \ldots) - \mathbf{x}| \leq \delta_2) \\
&\leq \lim_{m \to \infty} \lim_{\delta \to 0} \limsup_{\theta \to \infty} \frac{a(\theta)}{\theta} \\
&\qquad \times \log P\left(\max_{1 \leq i \leq m} |a(\theta)X_i(\theta) - x_i| \leq 2^m \delta\right).
\end{aligned} \tag{4.5}$$

It is known (cf. page 107 of [1]) that for $\sum_{i=1}^{n} y_i < 1$, where $0 \leq y_k < 1, k = 1, \ldots, n$, the joint density function of $(X_1(\theta), \ldots, X_n(\theta))$ is

$$\frac{\theta^n (1-(y_1+\cdots+y_n))^{\theta-1}}{(1-y_1)(1-(y_1+y_2))(1-(y_1+\cdots+y_{n-1}))}. \tag{4.6}$$

For any $n \geq 1, \delta > 0$, it follows from (4.6) that for sufficiently large $\theta$

$$\left\{\left(1 - \frac{x_1+\cdots+x_n+n\delta}{a(\theta)}\right)^{\theta-1} \left(\frac{\delta\theta}{a(\theta)}\right)^n\right\}$$
$$\times \left\{\left(1 - \frac{x_1-\delta}{a(\theta)}\right)\left(1 - \frac{x_1+x_2-2\delta}{a(\theta)}\right)\cdots\right.$$
$$\left.\times \left(1 - \frac{x_1+\cdots+x_{n-1}-(n-1)\delta}{a(\theta)}\right)\right\}^{-1}$$
$$\leq P\left(\max_{1\leq i\leq n} |a(\theta)X_i(\theta) - x_i| < \delta\right)$$
$$\leq \left\{\left(1 - \frac{x_1+\cdots+x_n-n\delta}{a(\theta)}\right)^{\theta-1} \left(\frac{\delta\theta}{a(\theta)}\right)^n\right\}$$
$$\times \left\{\left(1 - \frac{x_1+\delta}{a(\theta)}\right)\left(1 - \frac{x_1+x_2+2\delta}{a(\theta)}\right)\cdots\right.$$
$$\left.\times \left(1 - \frac{x_1+\cdots+x_{n-1}+(n-1)\delta}{a(\theta)}\right)\right\}^{-1}.$$

Therefore

$$\lim_{\delta\to 0} \limsup_{\theta\to\infty} \frac{a(\theta)}{\theta} \log P\left(\max_{1\leq i\leq n} |a(\theta)X_i(\theta) - x_i| \leq \delta\right)$$
$$= \lim_{\delta\to 0} \liminf_{\theta\to\infty} \frac{a(\theta)}{\theta} \log P\left(\max_{1\leq i\leq n} |a(\theta)X_i(\theta) - x_i| < \delta\right)$$
$$= -\sum_{i=1}^{n} x_i$$

which combined with (4.3) and (4.5) implies that

$$\lim_{\delta\to 0} \limsup_{\theta\to\infty} \frac{a(\theta)}{\theta} \log P(|a(\theta)(X_1(\theta), X_2(\theta), \ldots) - \mathbf{x}| \leq \delta)$$
$$= \lim_{\delta\to 0} \liminf_{\theta\to\infty} \frac{a(\theta)}{\theta} \log P(|a(\theta)(X_1(\theta), X_2(\theta), \ldots) - \mathbf{x}| < \delta)$$
$$= -\sum_{i=1}^{\infty} x_i.$$

Now we show the exponential tightness. For any $n \geq 1$ and $L \geq 1$, it follows from direct calculation that

$$P(a(\theta)X_n(\theta) \geq L) \leq P(a(\theta)U_n \geq L) = \left(1 - \frac{L}{a(\theta)}\right)_+^\theta,$$

where $(1 - \frac{L}{a(\theta)})_+$ is the positive part of $(1 - \frac{L}{a(\theta)})$. Set $K = \prod_{i=1}^\infty [0, iL]$. Then $K$ is a compact subset of $\mathbb{R}_+^\infty$. Noting that for $x \geq 0$

$$(1-x)_+ \leq e^{-x},$$

we get

$$P(a(\theta)(X_1(\theta), X_2(\theta), \ldots) \notin K) \leq \sum_{i=1}^\infty P(a(\theta)X_i(\theta) \geq iL)$$
$$\leq \sum_{i=1}^\infty \left(1 - \frac{iL}{a(\theta)}\right)_+^\theta$$
$$\leq \sum_{i=1}^\infty \exp\left\{-i\frac{\theta L}{a(\theta)}\right\}$$

which implies

$$\limsup_{\theta \to \infty} \frac{a(\theta)}{\theta} \log P(a(\theta)(X_1(\theta), X_2(\theta), \ldots) \notin K) \leq -L. \qquad \square$$

**5. Moderate deviations for the Poisson–Dirichlet distribution.** Theorem 3.3 says that $\mathbf{P}(\theta) = (P_1(\theta), P_2(\theta), \ldots)$ approaches a nontrivial random sequence when scaled by a factor of $\theta$ and shifted by $\beta(\theta)$. Replacing the scaling factor by $a(\theta)$ satisfying (3.5), we get

$$(5.1) \qquad a(\theta)\left(\mathbf{P}(\theta) - \frac{\beta(\theta)}{\theta}(1, 1, \ldots)\right) \to (0, 0, \ldots).$$

The LDP corresponds to the case when $a(\theta) = 1$ and has been established in [4]. In this section, we establish the MDP for $\mathbf{P}(\theta) = (P_1(\theta), P_2(\theta), \ldots)$ or, equivalently, the LDP associated with the limits in (5.1). Considering the connection to Poisson point process, it is thus natural to start with the MDP for $\frac{1}{\theta}(\sigma_1(\theta), \sigma_2(\theta), \ldots)$.

We first establish the MDP for $\sigma_n(\theta)/\theta$ for any $n$ followed by the MDP for $(\sigma_1(\theta), \ldots, \sigma_n(\theta))/\theta$. The infinite-dimensional case follows from finite-dimensional approximation. To go from the MDP for $\frac{1}{\theta}(\sigma_1(\theta), \sigma_2(\theta), \ldots)$ to the MDP for the Poisson–Dirichlet distribution, one would hope to prove that a certain exponential equivalency holds.

5.1. *MDP for $\frac{\sigma_n(\theta)}{\theta}$.* It is known (cf. [15]) that for each $n \geq 1$, the density function of $(\sigma_1(\theta), \ldots, \sigma_n(\theta))$ is

(5.2) $$f_n(u_1, \ldots, u_n) = \frac{\theta^n}{u_1 \cdots u_n} e^{-\sum_{i=1}^n u_i - \theta E_1(u_n)}, \qquad u_1 \geq u_2 \geq \cdots \geq u_n > 0.$$

In particular, the density function of $\sigma_1(\theta)$ is

(5.3) $$\theta u^{-1} e^{-u - \theta E_1(u)}, \qquad u > 0,$$

where $E_1(u) = \int_u^\infty y^{-1} e^{-y}\, dy$. We extend $E_1(u)$ to the whole real line by defining $E_1(u) = +\infty$ for $u \leq 0$.

The distribution function of $\sigma_1(\theta)$ is

(5.4) $$P\{\sigma_1(\theta) \leq u\} = e^{-\theta E_1(u)}, \qquad u > 0.$$

One can find on page 146 in [15] the following explicit expression for the distribution function of $\sigma_n(\theta)$ for all $n \geq 1$.

LEMMA 5.1. *The distribution function of $\sigma_n(\theta)$ is*

(5.5) $$F_n(y) = \frac{1}{(n-1)!} \int_{\theta E_1(y)}^\infty u^{n-1} e^{-u}\, du, \qquad y > 0.$$

Next we establish the MDP for $\sigma_1(\theta)/\theta$.

THEOREM 5.1. *The MDP holds for $\sigma_1(\theta)/\theta$ with speed $\frac{a(\theta)}{\theta}$ and rate function*

$$J_1(x) = \begin{cases} x, & x \geq 0, \\ \infty, & otherwise. \end{cases}$$

PROOF. For any fixed $x$, we have

(5.6) $$P\left\{ a(\theta) \left( \frac{\sigma_1(\theta) - \beta(\theta)}{\theta} \right) \leq x \right\} = e^{-\theta E_1((\theta/a(\theta))x + \beta(\theta))}.$$

By L'Hospital's rule,

(5.7) $$\lim_{x \to \infty} x e^x E_1(x) = \lim_{x \to \infty} \frac{x}{x+1} = 1.$$

Restricting to a subsequence if necessary we can assume without loss of generality that $\lim_{\theta \to \infty} [\frac{\theta}{a(\theta)} x + \beta(\theta))]$ exists in $[-\infty, +\infty]$. If the limit is negative, then the event $\{a(\theta)(\frac{\sigma_1(\theta) - \beta(\theta)}{\theta}) \leq x\}$ is eventually empty. Therefore

(5.8) $$\limsup_{\theta \to \infty} \frac{a(\theta)}{\theta} \log P\left( a(\theta)\left( \frac{\sigma_1(\theta) - \beta(\theta)}{\theta} \right) \leq x \right) = -\infty.$$

If $\lim_{\theta\to\infty}[\frac{\theta}{a(\theta)}x + \beta(\theta)]$ is a nonnegative finite number, then $x$ is negative and $\frac{\theta}{a(\theta)}$ and $\beta(\theta)$ are of the same scale as $\log \theta$. It follows from (5.6) that (5.8) also holds in this case.

When $\lim_{\theta\to\infty}[\frac{\theta}{a(\theta)}x + \beta(\theta)] = \infty$, we can use (5.7) to get

$$\lim_{\theta\to\infty} \frac{a(\theta)}{\theta}\left(-\theta E_1\left(\frac{\theta}{a(\theta)}x + \beta(\theta)\right)\right)$$

(5.9)
$$= -\lim_{\theta\to\infty} \frac{a(\theta)}{\theta} \frac{\log\theta}{(\theta/a(\theta))x + \beta(\theta)} e^{-(\theta/a(\theta))x}$$

$$= \begin{cases} 0, & x \geq 0 \\ -\infty, & x < 0. \end{cases}$$

Thus

(5.10) $\quad \lim_{\theta\to\infty} \frac{a(\theta)}{\theta} \log P\left(a(\theta)\left(\frac{\sigma_1(\theta) - \beta(\theta)}{\theta}\right) \leq x\right) = 0, \quad x \geq 0,$

and

(5.11) $\quad \limsup_{\theta\to\infty} \frac{a(\theta)}{\theta} \log P\left(a(\theta)\left(\frac{\sigma_1(\theta) - \beta(\theta)}{\theta}\right) \leq x\right) = -\infty, \quad x < 0.$

For $x \geq 0$, it follows from (5.6) and (5.7) that

$$\limsup_{\theta\to\infty} \frac{a(\theta)}{\theta} \log P\left(a(\theta)\left(\frac{\sigma_1(\theta) - \beta(\theta)}{\theta}\right) \geq x\right)$$

(5.12)
$$= \limsup_{\theta\to\infty} \frac{a(\theta)}{\theta} \log[1 - e^{-\theta E_1((a(\theta)/\theta)x + \beta(\theta))}]$$

$$= \limsup_{\theta\to\infty} \frac{a(\theta)}{\theta} \log\left[\theta E_1\left(\frac{a(\theta)}{\theta}x + \beta(\theta)\right)\right] \leq -x.$$

Together, (5.12) and (5.11) imply that the family of the laws of $a(\theta)(\frac{\sigma_1(\theta) - \beta(\theta)}{\theta})$ is exponentially tight.

Let $g_1(u)$ denote the density function of $a(\theta)(\frac{\sigma_1(\theta) - \beta(\theta)}{\theta})$. Then it follows from (5.3) that

(5.13) $\quad g_1(u) = \frac{\theta}{a(\theta)} \frac{\log\theta}{(\theta/a(\theta))u + \beta(\theta)} e^{-(\theta/a(\theta))u} e^{-\theta E_1((\theta/a(\theta))u + \beta(\theta))}.$

This, combined with (5.9), implies that

(5.14) $\quad\quad\quad\quad \frac{a(\theta)}{\theta} \log g_1(u) \to -u, \quad u > 0,$

(5.15) $\quad\quad\quad\quad \frac{a(\theta)}{\theta} \log g_1(u) \to -\infty, \quad u < 0.$

For each $x \neq 0$, choose $\delta$ small enough so that all numbers in the interval $[x-\delta, x+\delta]$ are of the same sign. It is not hard to see that for $u \in [x-\delta, x+\delta]$,

$$g_1(u) \geq \frac{\theta}{a(\theta)} \frac{\log \theta}{(\theta/a(\theta))(x+\delta) + \beta(\theta)}$$
(5.16)
$$\times e^{-(\theta/a(\theta))(x+\delta)} e^{-\theta E_1((\theta/a(\theta))(x-\delta)+\beta(\theta))}$$

and

$$g_1(u) \leq \frac{\theta}{a(\theta)} \frac{\log \theta}{(\theta/a(\theta))(x-\delta) + \beta(\theta)}$$
(5.17)
$$\times e^{-(\theta/a(\theta))(x-\delta)} e^{-\theta E_1((\theta/a(\theta))(x+\delta)+\beta(\theta))}.$$

Putting (5.14), (5.15), (5.16) and (5.17) together, we get that for $x > 0$,

$$\lim_{\delta \to 0} \limsup_{\theta \to \infty} \frac{a(\theta)}{\theta} \log P\left(a(\theta)\left(\frac{\sigma_1(\theta) - \beta(\theta)}{\theta}\right) \in (x - \delta, x + \delta)\right)$$
$$= \lim_{\delta \to 0} \liminf_{\theta \to \infty} \frac{a(\theta)}{\theta} \log P\left(a(\theta)\left(\frac{\sigma_1(\theta) - \beta(\theta)}{\theta}\right) \in (x - \delta, x + \delta)\right) = -x,$$

and for any $x < 0$,

$$\lim_{\delta \to 0} \limsup_{\theta \to \infty} \frac{a(\theta)}{\theta} \log P\left(a(\theta)\left(\frac{\sigma_1(\theta) - \beta(\theta)}{\theta}\right) \in (x - \delta, x + \delta)\right) = -\infty.$$

Together, (5.10) and (5.11) imply that

$$\lim_{\delta \to 0} \limsup_{\theta \to \infty} \frac{a(\theta)}{\theta} \log P\left(a(\theta)\left(\frac{\sigma_1(\theta) - \beta(\theta)}{\theta}\right) \in (-\delta, \delta)\right)$$
$$= \lim_{\delta \to 0} \liminf_{\theta \to \infty} \frac{a(\theta)}{\theta} \log P\left(a(\theta)\left(\frac{\sigma_1(\theta) - \beta(\theta)}{\theta}\right) \in (-\delta, \delta)\right) = 0.$$

The theorem now follows from the local LDP and exponential tightness. □

The next theorem gives the MDP of $\sigma_n(\theta)/\theta$ for $n \geq 2$.

THEOREM 5.2. *The MDP holds for $\sigma_n(\theta)/\theta$ with speed $\frac{a(\theta)}{\theta}$ and rate function $J_n(x) = nx, x \geq 0$.*

PROOF. For $x > 0$, it follows from (5.10) that

$$0 \geq \lim_{\theta \to \infty} \frac{a(\theta)}{\theta} \log P\left(a(\theta)\left(\frac{\sigma_n(\theta) - \beta(\theta)}{\theta}\right) \leq x\right)$$
(5.18)
$$\geq \lim_{\theta \to \infty} \frac{a(\theta)}{\theta} \log P\left(a(\theta)\left(\frac{\sigma_1(\theta) - \beta(\theta)}{\theta}\right) \leq x\right) = 0.$$

By L'Hospital's rule,

(5.19) $$\lim_{y \to 0} n y^{-n} e^y \int_0^y u^{n-1} e^{-u} \, du = 1.$$

Thus it follows from Lemma 5.1, (5.9) and (5.19) that

(5.20)
$$\begin{aligned}
\lim_{\theta \to \infty} \frac{a(\theta)}{\theta} \log P\bigg(a(\theta)\bigg(\frac{\sigma_n(\theta) - \beta(\theta)}{\theta}\bigg) \geq x\bigg) \\
= \lim_{\theta \to \infty} \frac{a(\theta)}{\theta} \log\bigg(1 - F_n\bigg(\frac{\theta}{a(\theta)} x + \beta(\theta)\bigg)\bigg) \\
= \lim_{\theta \to \infty} \frac{a(\theta)}{\theta} \log e^{-\theta E_1((\theta/a(\theta))x + \beta(\theta))} \\
+ \lim_{\theta \to \infty} \frac{a(\theta)}{\theta} n \log\bigg(\theta E_1\bigg(\frac{\theta}{a(\theta)} x + \beta(\theta)\bigg)\bigg) \\
= n \lim_{\theta \to \infty} \frac{a(\theta)}{\theta} \log\bigg(\frac{\log \theta}{(\theta/a(\theta))x + \beta(\theta)} e^{-(\theta/a(\theta))x}\bigg) \\
= -nx.
\end{aligned}$$

For $x < 0$, as in the proof of Theorem 5.1, it suffices to obtain estimates for those $x$ such that

$$\lim_{\theta \to \infty} \bigg(\frac{\theta}{a(\theta)} x + \beta(\theta)\bigg) = +\infty.$$

Since $\theta E_1(\frac{\theta}{a(\theta)} x + \beta(\theta)) \approx \frac{\log \theta}{(\theta/a(\theta))x + \beta(\theta)} e^{-(\theta/a(\theta))x}$ approaches infinity as $\theta$ tends to infinity, one gets that

(5.21)
$$\begin{aligned}
\limsup_{\theta \to \infty} \frac{a(\theta)}{\theta} \log P\bigg(a(\theta)\bigg(\frac{\sigma_n(\theta) - \beta(\theta)}{\theta}\bigg) \leq x\bigg) \\
= \limsup_{\theta \to \infty} \frac{a(\theta)}{\theta} \log\bigg(e^{-\theta E_1((\theta/a(\theta))x + \beta(\theta))} \\
\times \bigg(\theta E_1\bigg(\frac{\theta}{a(\theta)} x + \beta(\theta)\bigg)\bigg)^n\bigg) \\
= -\infty.
\end{aligned}$$

The exponential tightness of the laws of $\{a(\theta) \frac{\sigma_n(\theta) - \beta(\theta)}{\theta}\}$ now follows from (5.20) and (5.21). The local LDP can be obtained by an argument similar to that used in Theorem 5.1. $\square$

5.2. *MDP for* $\frac{1}{\theta}(\sigma_1(\theta), \sigma_2(\theta), \ldots)$. For each $n \geq 2$, we have

THEOREM 5.3. *The family* $\{P(a(\theta)(\frac{\sigma_1(\theta)-\beta(\theta)}{\theta},\ldots,\frac{\sigma_n(\theta)-\beta(\theta)}{\theta}) \in \cdot):\theta > 0\}$ *satisfies a LDP on* $\mathbb{R}^n$ *with speed* $\frac{a(\theta)}{\theta}$ *and rate function*

$$I_n(x_1,\ldots,x_n) = \begin{cases} \sum_{i=1}^n x_i, & \text{if } 0 \leq x_n \leq \cdots \leq x_1, \\ +\infty, & \text{otherwise.} \end{cases} \tag{5.22}$$

PROOF. It follows from (5.2) that for $x_1 \geq x_2 \cdots \geq x_n$ and $\frac{\theta}{a(\theta)}x_n + \beta(\theta) > 0$, the density function of $\frac{a(\theta)}{\theta}(\sigma_1(\theta) - \beta(\theta),\ldots,\sigma_n(\theta) - \beta(\theta))$ is

$$g_n(x_1,\ldots,x_n) = \left(\frac{\theta}{a(\theta)}\right)^n \left(\prod_{i=1}^n \frac{\log(\theta)}{(\theta/a(\theta))x_i + \beta(\theta)}\right) \\ \times e^{-[(\theta/a(\theta))\sum_{i=1}^n x_i + \theta E_1((\theta/a(\theta))x_n + \beta(\theta))]}. \tag{5.23}$$

By direct calculation,

$$\frac{a(\theta)}{\theta} \log g_n(x_1,\ldots,x_n) \to -\sum_{i=1}^n x_i, \qquad x_n > 0, \tag{5.24}$$

$$\frac{a(\theta)}{\theta} \log g_n(x_1,\ldots,x_n) \to -\infty, \qquad x_n < 0. \tag{5.25}$$

For $x_1 \geq x_2 \cdots \geq x_n$, let $B((x_1,\ldots,x_n),\delta)$ denote the closed ball centered at $(x_1,\ldots,x_n)$ with radius $\delta$, and $B^\circ((x_1,\ldots,x_n),\delta)$ be the corresponding open ball. Then for $x_n > 0$,

$$\lim_{\delta \to 0} \limsup_{\theta \to \infty} \frac{a(\theta)}{\theta} \log P\Big(\frac{a(\theta)}{\theta}(\sigma_1(\theta) - \beta(\theta),\ldots,\sigma_n(\theta) - \beta(\theta)) \\ \in B((x_1,\ldots,x_n),\delta)\Big)$$

$$= \lim_{\delta \to 0} \liminf_{\theta \to \infty} \frac{a(\theta)}{\theta} \log P\Big(\frac{a(\theta)}{\theta}(\sigma_1(\theta) - \beta(\theta),\ldots,\sigma_n(\theta) - \beta(\theta)) \\ \in B^\circ((x_1,\ldots,x_n),\delta)\Big) \tag{5.26}$$

$$= -\sum_{i=1}^n x_i,$$

and for any $x_n < 0$,

$$\lim_{\delta \to 0} \limsup_{\theta \to \infty} \frac{a(\theta)}{\theta} \log P\Big(\frac{a(\theta)}{\theta}(\sigma_1(\theta) - \beta(\theta),\ldots,\sigma_n(\theta) - \beta(\theta)) \\ \in B((x_1,\ldots,x_n),\delta)\Big)$$

$$(5.27) \qquad = \lim_{\delta \to 0} \liminf_{\theta \to \infty} \frac{a(\theta)}{\theta} \log P\Big(\frac{a(\theta)}{\theta}(\sigma_1(\theta) - \beta(\theta)), \ldots, \sigma_n(\theta) - \beta(\theta))$$

$$\in B^\circ((x_1, \ldots, x_n), \delta)\Big)$$

$$= -\infty.$$

If $x_1 = 0$, the upper estimate follows from Theorem 5.1. If $x_{r-1} > 0, x_r = 0$ for some $1 < r \leq n$, then the upper estimate is obtained from that of $\frac{a(\theta)}{\theta}(\sigma_1(\theta) - \beta(\theta), \ldots, \sigma_{r-1}(\theta) - \beta(\theta))$. The lower estimate when $x_r = 0$ for some $1 \leq r \leq n$ is obtained by approximating the boundary with open subsets that have all positive coordinates.

Fix an $L > 0$. Noting that $\bigcup_{i=1}^n \{\frac{a(\theta)}{\theta}(\sigma_i(\theta) - \beta(\theta)) > L\} = \{\frac{a(\theta)}{\theta}(\sigma_1(\theta) - \beta(\theta)) > L\}$, it follows that

$$(5.28) \qquad \lim_{L \to \infty} \limsup_{\theta \to \infty} \frac{a(\theta)}{\theta} \log P\left\{\bigcup_{i=1}^n \left\{\frac{a(\theta)}{\theta}(\sigma_i(\theta) - \beta(\theta)) > L\right\}\right\} = -\infty.$$

On the other hand,

$$(5.29) \qquad \begin{aligned} &\limsup_{\theta \to \infty} \frac{a(\theta)}{\theta} \log P\left\{\bigcup_{i=1}^n \left\{\frac{a(\theta)}{\theta}(\sigma_i(\theta) - \beta(\theta)) < -L\right\}\right\} \\ &\leq \limsup_{\theta \to \infty} \frac{a(\theta)}{\theta} \log P\left\{\frac{a(\theta)}{\theta}(\sigma_n(\theta) - \beta(\theta)) \leq -L\right\} = -\infty. \end{aligned}$$

Therefore we have the exponential tightness and the theorem. $\square$

The MDP for $\frac{1}{\theta}(\sigma_1(\theta), \sigma_2(\theta), \ldots)$ is derived in the next theorem.

THEOREM 5.4. *The family $\{P(\frac{a(\theta)}{\theta}(\sigma_1(\theta) - \beta(\theta), \sigma_2(\theta) - \beta(\theta), \ldots) \in \cdot) : \theta > 0\}$ satisfies a LDP on $\mathbb{R}^\infty$ with speed $\frac{a(\theta)}{\theta}$ and rate function*

$$(5.30) \qquad I(x_1, x_2, \ldots) = \begin{cases} \sum_{i=1}^\infty x_i, & x_1 \geq \cdots \geq 0, \\ \infty, & \textit{otherwise.} \end{cases}$$

PROOF. Identify $\mathbb{R}^\infty$ with the projective limit of $\mathbb{R}^n, n = 1, \ldots$. Then the theorem follows from Theorem 3.3 in [5] and Theorem 5.3. $\square$

5.3. *MDP for the Poisson–Dirichlet distribution.* Using the results in the previous subsection we now derive the MDP for the Poisson–Dirichlet distribution. The representation (1.1), combined with the fact that $\sigma(\theta)$ is approximately $\theta$, seems to suggest that the MDP for the Poisson–Dirichlet distribution should follow from the MDP for $\frac{1}{\theta}(\sigma_1(\theta), \sigma_2(\theta), \ldots)$. This turns out to be true. It does not seem to be easy to get a more direct proof using the explicit expression in [28] of the density functions of $(P_1(\theta), \ldots, P_n(\theta))$ for each $n \geq 1$.

THEOREM 5.5. *For each $n \geq 1$, the family $\{P(a(\theta)(P_1(\theta) - \frac{\beta(\theta)}{\theta}, \ldots, P_n(\theta) - \frac{\beta(\theta)}{\theta}, \ldots) \in \cdot) : \theta > 0\}$ satisfies a LDP on $\mathbb{R}^\infty$ with speed $\frac{a(\theta)}{\theta}$ and rate function*

$$(5.31) \qquad I(x_1, x_2, \ldots) = \begin{cases} \sum_{i=1}^{\infty} x_i, & x_1 \geq \cdots \geq 0, \\ \infty, & \text{otherwise.} \end{cases}$$

PROOF. From representation (1.1), one obtains that

$$(5.32) \qquad \begin{aligned} & a(\theta)\left(P_n(\theta) - \frac{\beta(\theta)}{\theta}\right) \\ &= \frac{\theta}{\sigma(\theta)} a(\theta) \left[\frac{\sigma_n(\theta) - \beta(\theta)}{\theta}\right] + \frac{a(\theta)\beta(\theta)}{\theta}\left(\frac{\theta}{\sigma(\theta)} - 1\right). \end{aligned}$$

Write

$$\gamma(\theta) = \frac{a(\theta)\beta(\theta)}{\theta},$$

and without loss of generality we assume that

$$\lim_{\theta \to \infty} \gamma(\theta) = c \in [0, +\infty].$$

It is clear that

$$(5.33) \qquad \frac{a(\theta)}{\gamma^2(\theta)} = \frac{\theta^2}{a(\theta)\beta^2(\theta)} \to \infty, \qquad \theta \to \infty.$$

If $c < \infty$, it follows from Corollary 3.1 that for any $L > 0$

$$(5.34) \qquad \limsup_{\theta \to \infty} \frac{a(\theta)}{\theta} \log P\left\{\gamma(\theta)\left|\frac{\theta}{\sigma(\theta)} - 1\right| \geq L\right\} = -\infty.$$

For $c = \infty$, and any $1 > \delta > 0$

$$(5.35) \qquad \begin{aligned} & \left\{\gamma(\theta)\left|\frac{\theta}{\sigma(\theta)} - 1\right| \geq L\right\} \\ & \subset \left\{\gamma(\theta)\left|\frac{\sigma(\theta)}{\theta} - 1\right| \geq L(1-\delta)\right\} \cup \left\{\left|\frac{\sigma(\theta)}{\theta} - 1\right| \geq \delta\right\}. \end{aligned}$$

Since $\gamma(\theta) \leq \beta(\theta)$ and $\lim_{\theta \to \infty} \frac{\beta(\theta)}{\sqrt{\theta}} = 0$, it follows from the MDP (Theorem 3.2) for $\sigma(\theta)/\theta$, and (5.33) that

(5.36)
$$\limsup_{\theta \to \infty} \frac{a(\theta)}{\theta} \log P\left\{\gamma(\theta) \left|\frac{\sigma(\theta)}{\theta} - 1\right| \geq (1-\delta)L\right\}$$
$$= \limsup_{\theta \to \to \infty} \frac{a(\theta)}{\gamma^2(\theta)} \frac{\gamma^2(\theta)}{\theta} \log P\left\{\gamma(\theta)\left|\frac{\sigma(\theta)}{\theta} - 1\right| \geq (1-\delta)L\right\} = -\infty,$$

which, combined with Corollary 3.1 and (5.35), shows that (5.34) still holds in this case. Therefore $a(\theta)(P_n(\theta) - \frac{\beta(\theta)}{\theta})$ and $\frac{\theta}{\sigma(\theta)}a(\theta)[\frac{\sigma_n(\theta)-\beta(\theta)}{\theta}]$ are exponentially equivalent.

Since $\frac{\theta}{\sigma(\theta)}a(\theta)[\frac{\sigma_n(\theta)-\beta(\theta)}{\theta}]$ is exponentially equivalent to $a(\theta)[\frac{\sigma_n(\theta)-\beta(\theta)}{\theta}]$ by Lemma 2.1 and Corollary 3.1, it follows that $a(\theta)(P_n(\theta) - \frac{\beta(\theta)}{\theta})$ and $a(\theta)[\frac{\sigma_n(\theta)-\beta(\theta)}{\theta}]$ are exponentially equivalent for all $n \geq 1$. Thus the MDP for the Poisson–Dirichlet distribution is the same as the MDP for $\frac{1}{\theta}(\sigma_1(\theta), \sigma_2(\theta), \ldots)$. □

**6. Moderate deviations for homozygosity.** For each $m \geq 2$, it is shown in Theorem 3.4 that the scaled homozygosity $\frac{\theta^{m-1}}{\Gamma(m)}H_m(\mathbf{P}(\theta))$ satisfies a fluctuation theorem with $c(\theta) = 1$ and $b(\theta) = \sqrt{\theta}$. It is thus natural to consider MDPs for $\frac{\theta^{m-1}}{\Gamma(m)}H_m(\mathbf{P}(\theta))$ or equivalently the LDP for the family $\{a(\theta)[\frac{\theta^{m-1}}{\Gamma(m)}H_m(\mathbf{P}(\theta)) - 1] : \theta > 0\}$ for a scale $a(\theta)$ satisfying

(6.1) $$\lim_{\theta \to \infty} a(\theta) = \infty, \qquad \lim_{\theta \to \infty} \frac{a(\theta)}{\sqrt{\theta}} = 0.$$

It will turn out in Remark 6.1 that the following additional restriction on $a(\theta)$ is necessary in order to get the MDP: for some $0 < \varepsilon < 1/(2m-1)$,

(6.2) $$\liminf_{\theta \to \infty} \frac{a^{1-\varepsilon}(\theta)}{\theta^{(m-1)/(2m-1)}} > 0.$$

The main idea of the proof is to explore the connection between homozygosity and the Poisson process, and apply Campbell's theorem.

Let us first consider the MDP of

$$G_\theta^{(m)} := \sum_{j=1}^\infty \sigma_j^m(\theta).$$

It follows from Campbell's theorem that

$$E(e^{it \sum_{j=1}^\infty \sigma_j^m(\theta)}) = \exp\left\{\theta \int_0^\infty (e^{ity^m} - 1)y^{-1}e^{-y}\,dy\right\}$$

which implies that $\{G_\theta^{(m)}, \theta \geq 0\}$ is a random process with stationary and independent increments. The difficulty here is that the exponential moment is not finite. MDPs for models with infinite exponential moment have been studied in [14, 17, 22]. A typical way of establishing the MDP is to verify the following Ledoux condition [22]: there exists a constant $M > 0$ such that for any $\delta > 0$,

$$\limsup_{\theta \to \infty} \frac{a^2(\theta)}{\theta} \log\left(\theta P\left(|G_1^{(m)} - E(G_1^{(m)})| > \frac{\delta \theta}{a(\theta)}\right)\right) \leq -\frac{\delta^2}{M}.$$

This condition does not seem to be easy to verify for our model. Therefore we employ a truncation procedure.

LEMMA 6.1. *Set*

$$G_\theta = (\sigma(\theta) - \theta, G_\theta^{(m)} - \Gamma(m)\theta).$$

*Then the family* $\{\frac{a(\theta)}{\theta} G_\theta : \theta > 0\}$ *satisfies a LDP with speed* $\frac{a^2(\theta)}{\theta}$ *and rate function*

$$\Lambda^*(x, y) := \frac{1}{2(\Gamma(2m) - \Gamma(m+1)^2)}(\Gamma(2m)x^2 - 2\Gamma(m+1)xy + y^2),$$

$$x \in \mathbb{R}, y \in \mathbb{R}.$$

PROOF. By (6.1) and (6.2), there exist $\tau > 0$ and a positive integer $l \geq 3 \vee \frac{2}{(2m-1)\varepsilon}$ such that

$$\lim_{\theta \to \infty} \frac{a(\theta)}{\theta^\tau} = +\infty$$

and

$$\lim_{\theta \to \infty} \left(\frac{a^2(\theta)}{\theta} a^{(l-2)/(m-1)l}(\theta)\right)^{(m-1)l} = \lim_{\theta \to \infty} \frac{a^{(2m-1)l-2}(\theta)}{\theta^{(m-1)l}}$$

$$= \lim_{\theta \to \infty} \frac{a^{1-2/((2m-1)l)}(\theta)}{\theta^{(m-1)/(2m-1)}} = \infty.$$

Take

$$\gamma(\theta) = \frac{a^{(l-2)/(m-1)l}(\theta)}{\log((a^2(\theta)/\theta)a^{(l-2)/(m-1)l}(\theta))}.$$

Then $\gamma(\theta)$ grows faster than a positive power of $\theta$ and

$$\lim_{\theta \to \infty} \frac{\gamma(\theta)}{a^{(l-2)/(m-1)l}(\theta)} = 0, \quad \lim_{\theta \to \infty} \frac{a^2(\theta)\gamma(\theta)}{\theta} = \infty.$$

Set
$$\tilde{G}_\theta^{(1)} = \sum_{j=1}^\infty \sigma_j(\theta) I_{\{\sigma_j(\theta) \leq \gamma(\theta)\}}, \qquad \tilde{G}_\theta^{(m)} = \sum_{j=1}^\infty \sigma_j^m(\theta) I_{\{\sigma_j(\theta) \leq \gamma(\theta)\}}$$

and
$$\tilde{G}_\theta = (\tilde{G}_\theta^{(1)} - E(\tilde{G}_\theta^{(1)}), \tilde{G}_\theta^{(m)} - E(\tilde{G}_\theta^{(m)})).$$

Define
$$\Lambda(\alpha, \beta) = \tfrac{1}{2}(\alpha^2 + 2\Gamma(m+1)\alpha\beta + \Gamma(2m)\beta^2)$$
$$= \tfrac{1}{2} \begin{pmatrix} \alpha & \beta \end{pmatrix} \begin{pmatrix} 1 & \Gamma(m+1) \\ \Gamma(m+1) & \Gamma(2m) \end{pmatrix} \begin{pmatrix} \alpha \\ \beta \end{pmatrix}, \qquad \alpha \in \mathbb{R}, \beta \in \mathbb{R}.$$

Then
$$\sup_{\alpha \in \mathbb{R}, \beta \in \mathbb{R}} \{\alpha x + \beta y - \Lambda(\alpha, \beta)\}$$
$$= \frac{1}{2(\Gamma(2m) - \Gamma(m+1)^2)} \begin{pmatrix} x & y \end{pmatrix} \begin{pmatrix} \Gamma(2m) & -\Gamma(m+1) \\ -\Gamma(m+1) & 1 \end{pmatrix} \begin{pmatrix} x \\ y \end{pmatrix}$$
$$= \frac{1}{2(\Gamma(2m) - \Gamma(m+1)^2)} (\Gamma(2m)x^2 - 2\Gamma(m+1)xy + y^2),$$
$$x \in \mathbb{R}, y \in \mathbb{R}.$$

For any $\alpha \in \mathbb{R}$, $\beta \in \mathbb{R}$,
$$\left| \frac{a^2(\theta)}{\theta} \log E\left( \exp\left\{ \frac{1}{a(\theta)} (\alpha(\tilde{G}_\theta^{(1)} - E(\tilde{G}_\theta^{(1)})) + \beta(\tilde{G}_\theta^{(m)} - E(\tilde{G}_\theta^{(m)}))) \right\} \right) \right.$$
$$\left. - \Lambda(\alpha, \beta) \right|$$
$$= \left| \frac{a^2(\theta)}{\theta} \log \exp\left\{ \theta \int_0^{\gamma(\theta)} \left( e^{(1/a(\theta))(\alpha y + \beta y^m)} \right. \right. \right.$$
$$\left. \left. \left. - 1 - \frac{1}{a(\theta)}(\alpha y + \beta y^m) \right) y^{-1} e^{-y} \, dy \right\} - \Lambda(\alpha, \beta) \right|$$
$$= \left| a^2(\theta) \int_0^{\gamma(\theta)} (e^{(1/a(\theta))(\alpha y + \beta y^m)} - 1 - a^{-1}(\theta)(\alpha y + \beta y^m)) y^{-1} e^{-y} \, dy \right.$$
$$\left. - \Lambda(\alpha, \beta) \right|$$
$$\leq \left| \int_0^{\gamma(\theta)} \frac{1}{2} (\alpha y + \beta y^m)^2 y^{-1} e^{-y} \, dy - \Lambda(\alpha, \beta) \right|$$

$$+ \sum_{k=3}^{l} \frac{1}{k!} a^{-(k-2)}(\theta) \int_0^{\gamma(\theta)} |\alpha y + \beta y^m|^k y^{-1} e^{-y}\, dy$$

$$+ \sum_{k=l+1}^{\infty} \frac{1}{k!} a^{-(k-2)}(\theta)(|\alpha| + |\beta|\gamma(\theta)^{m-1})^k \Gamma(k)$$

$$\to 0 \quad \text{as } \theta \to \infty.$$

Therefore, by the Gärtner–Ellis theorem, $(\frac{a(\theta)}{\theta}\tilde{G}_\theta, \frac{a^2(\theta)}{\theta}, \Lambda^*)$ satisfies a LDP. Noting that $\gamma(\theta)$ grows faster than $\theta^\alpha$ for a certain $\alpha > 0$, it follows that

$$\lim_{\theta \to \infty} \theta E_1(\gamma(\theta)) = \lim_{\theta \to \infty} \theta \gamma(\theta)^{-1} e^{-\gamma(\theta)} = 0.$$

Taking into account the fact that $\tilde{G}_\theta^{(1)} \leq \sigma(\theta)$, $\tilde{G}_\theta^{(m)} \leq G_\theta^{(m)}$, we have that

$$\limsup_{\theta \to \infty} \frac{a^2(\theta)}{\theta} \log P\Big(|(\sigma(\theta), G_\theta^{(m)}) - (\tilde{G}_\theta^{(1)}, \tilde{G}_\theta^{(m)})| \geq \delta \frac{\theta}{a(\theta)}\Big)$$

$$\leq \limsup_{\theta \to \infty} \frac{a^2(\theta)}{\theta} \log P\Big(|(\sigma(\theta), G_\theta^{(m)})| I_{\{\sigma_1(\theta) \geq \gamma(\theta)\}} \geq \delta \frac{\theta}{a(\theta)}\Big)$$

$$\leq \limsup_{\theta \to \infty} \frac{a^2(\theta)}{\theta} \log P(\sigma_1(\theta) \geq \gamma(\theta))$$

$$= \limsup_{\theta \to \infty} \frac{a^2(\theta)}{\theta} \log(1 - e^{-\theta E_1(\gamma(\theta))})$$

$$= \limsup_{\theta \to \infty} \frac{[\log \theta + \log E_1(\gamma(\theta))] a^2(\theta)}{\theta}$$

$$\leq \limsup_{\theta \to \infty} \frac{(\log \theta - \gamma(\theta)) a^2(\theta)}{\theta}$$

$$= -\infty$$

which implies that for any $\delta > 0$,

$$\limsup_{\theta \to \infty} \frac{a^2(\theta)}{\theta} \log P\Big(\frac{|G_\theta - \tilde{G}_\theta - (E(\tilde{G}_\theta^{(1)}) - \theta, E(\tilde{G}_\theta^{(m)}) - \Gamma(m)\theta)|}{\theta/a(\theta)} \geq \delta\Big)$$
(6.3)
$$= -\infty.$$

By direct calculation,

$$\lim_{\theta \to \infty} \frac{(E(\tilde{G}_\theta^{(1)}), E(\tilde{G}_\theta^{(m)})) - (\theta, \Gamma(m)\theta)}{\theta/a(\theta)}$$

$$= -\lim_{\theta \to \infty} a(\theta)\Big(\int_{\gamma(\theta)}^\infty e^{-y}\, dy, \int_{\gamma(\theta)}^\infty y^{2m-1} e^{-y}\, dy\Big)$$

$$= -\lim_{\theta \to \infty} a(\theta)(e^{-\gamma(\theta)}, \gamma^{2m-1}(\theta)e^{-\gamma(\theta)}) = 0,$$

which, combined with (6.3), implies that $\frac{a(\theta)}{\theta}\tilde{G}_\theta$ and $\frac{a(\theta)}{\theta}G_\theta$ are exponentially equivalent. Therefore

$$\left(\frac{a(\theta)}{\theta}G_\theta, \frac{a^2(\theta)}{\theta}, \Lambda^*\right)$$

satisfies the LDP. □

Now we are ready to prove the main result of this section.

THEOREM 6.1. *The family* $a(\theta)(\frac{\theta^{m-1}}{\Gamma(m)}H_m(\mathbf{P}(\theta))-1)$ *satisfies a LDP with speed* $\frac{a^2(\theta)}{\theta}$ *and rate function* $\frac{z^2}{2(\Gamma(2m)/\Gamma(m)^2-m^2)}$.

PROOF. By direct calculation,

$$a(\theta)\left(\frac{\theta^{m-1}}{\Gamma(m)}H_m(\mathbf{P}(\theta))-1\right)$$

$$= a(\theta)\left(\frac{\theta^{m-1}G_\theta^{(m)}}{\sigma^m(\theta)\Gamma(m)}-1\right)$$

$$= a(\theta)\left(\left(\frac{\theta}{\sigma(\theta)}\right)^m-1\right) + \left(\frac{\theta}{\sigma(\theta)}\right)^m\frac{a(\theta)(G_\theta^{(m)}-E(G_\theta^{(m)}))}{\Gamma(m)\theta}$$

$$= \frac{a(\theta)}{\theta}(\theta-\sigma(\theta))\sum_{k=1}^{m}\left(\frac{\theta}{\sigma(\theta)}\right)^k + \left(\frac{\theta}{\sigma(\theta)}\right)^m\frac{a(\theta)(G_\theta^{(m)}-E(G_\theta^{(m)}))}{\Gamma(m)\theta}.$$

Noting that for any $i \geq 1$ and for any $\delta > 0$,

$$\lim_{\theta \to \infty}\frac{a^2(\theta)}{\theta}\log P\left(\left|\left(\frac{\theta}{\sigma(\theta)}\right)^i-1\right| \geq \delta\right) = -\infty.$$

It then follows that

$$a(\theta)\left(\left(\frac{\theta}{\sigma(\theta)}\right)^m-1\right) \quad \text{and} \quad \left(\frac{\theta}{\sigma(\theta)}\right)^m\frac{a(\theta)(G_\theta^{(m)}-E(G_\theta^{(m)}))}{\Gamma(m)\theta},$$

are exponentially equivalent to

$$\frac{a(\theta)m(\theta-\sigma(\theta))}{\theta} \quad \text{and} \quad \frac{a(\theta)(G_\theta^{(m)}-E(G_\theta^{(m)}))}{\Gamma(m)\theta},$$

respectively. Thus

$$a(\theta)\left(\frac{\theta^{m-1}}{\Gamma(m)}H_m(\mathbf{P}(\theta))-1\right)$$

and

$$\frac{a(\theta)m(\theta - \sigma(\theta))}{\theta} + \frac{a(\theta)(G_\theta^{(m)} - E(G_\theta^{(m)}))}{\Gamma(m)\theta}$$

have the same LDP.

Since

$$\inf_{(y/\Gamma(m))-mx=z} \Lambda^*(x,y) = \frac{z^2}{2(\Gamma(2m)/\Gamma(m)^2 - m^2)},$$

Lemma 6.1 and the contraction principle yield that

$$\left(\frac{a(\theta)m(\theta - \sigma(\theta))}{\theta} + \frac{a(\theta)(G_\theta^{(m)} - E(G_\theta^{(m)}))}{\Gamma(m)\theta}, \frac{a^2(\theta)}{\theta}, \frac{z^2}{2(\Gamma(2m)/\Gamma(m)^2 - m^2)}\right)$$

satisfies a LDP, and the theorem follows. $\square$

REMARK 6.1. Choose the scaling factor, $a(\theta) = \theta^\gamma$. Then the MDP obtained here requires that $\gamma$ lies between $\frac{m-1}{2m-1}$ and $\frac{1}{2}$. It is natural to ask what happens for $\gamma \leq \frac{m-1}{2m-1}$. It follows from Lemma 6.1 and the contraction principle that the family $\{\frac{a(\theta)(G_\theta^{(m)} - E(G_\theta^{(m)}))}{\theta} : \theta > 0\}$ satisfies a LDP with speed $\frac{a^2(\theta)}{\theta}$ and a rate function $J(x) = \frac{x^2}{2\Gamma(2m)}$. Thus for any $\delta > 0$, there exists $\theta_0 > 0$ such that for all $\theta \geq \theta_0$,

$$P\left(|G_\theta^{(m)} - E(G_\theta^{(m)})| > \frac{\delta\theta}{a(\theta)}\right) \leq \exp\left\{-\frac{\theta}{a^2(\theta)}[J(\delta) - 1/2]\right\}.$$

Since $\{G_\theta^{(m)}, \theta \geq 0\}$ is a random process with stationary and independent increments, one can find sufficiently small $\delta_1, \delta_2 > 0$ such that

$$P\left(|G_1^{(m)} - E(G_1^{(m)})| > \frac{\delta\theta}{a(\theta)}\right)$$

$$\leq P\left(|G_{\theta+1}^{(m)} - E(G_{\theta+1}^{(m)})| > \frac{\delta_1(\theta+1)}{a(\theta+1)}\right) + P\left(|(G_\theta^{(m)} - E(G_\theta^{(m)}))| > \frac{\delta_2\theta}{a(\theta)}\right)$$

$$\leq 2\exp\left\{-\frac{\theta}{a^2(\theta)}[J(\delta_1 \wedge \delta_2) - 1/2]\right\}.$$

The fact that $\lim_{x \to \infty} J(x) = +\infty$, yields

$$\limsup_{\delta \to \infty} \limsup_{\theta \to \infty} \frac{a^2(\theta)}{\theta} \log P\left(|G_1^{(m)} - E(G_1^{(m)})| > \frac{\delta\theta}{a(\theta)}\right) = -\infty$$

which, combined with the fact that $E(G_1^{(m)})$ is a finite number, implies

(6.4) $$\limsup_{\delta \to \infty} \limsup_{\theta \to \infty} \frac{a^2(\theta)}{\theta} \log P\left(G_1^{(m)} > \frac{\delta\theta}{a(\theta)}\right) = -\infty.$$

Since

$$(6.5) \quad P\Big(G_1^{(m)} \geq \frac{\delta\theta}{a(\theta)}\Big) \geq P\Big(\sigma_1(1) \geq \Big(\frac{\delta\theta}{a(\theta)}\Big)^{1/m}\Big) = 1 - e^{-E_1((\delta\theta/a(\theta))^{1/m})},$$

it follows from (6.4) that

$$\limsup_{\delta\to\infty}\limsup_{\theta\to\infty} \frac{a^2(\theta)}{\theta} \log E_1\Big(\Big(\frac{\delta\theta}{a(\theta)}\Big)^{1/m}\Big) = -\infty.$$

Using the relation (5.7), one gets

$$\gamma > \frac{m-1}{2m-1},$$

which corresponds to the critical case of $\varepsilon = 0$ in (6.2). Thus the range of scaling obtained here is the best that one can get for the MDP with speed $\frac{a^2(\theta)}{\theta}$.

**7. Moderate deviations for the Dirichlet process.** In this section, the MDP for the Dirichlet process is derived through a combination of the LDP for the gamma distribution and MDPs for processes with stationary independent increments.

The Dirichlet$(\theta, \nu)$ distribution can be represented by

$$Z_{\theta,\nu}(t) = \frac{X(\theta\nu([0,t]))}{X(\theta)},$$

where $\{X(t), t \in [0,\infty)\}$ is a Gamma process. By Theorem 3.1, the family $\{P(X(\theta)/\theta \in \cdot), \theta > 0\}$ satisfies a LDP in $\mathbb{R}_+$ with speed $\theta$ and rate function $I(x)$ given by (3.4).

Let $a(\theta)$ be a positive function satisfying (6.1). With a time deformation, the following theorem is a minor generalization of the result in [25]. For completeness, a sketched proof is included.

THEOREM 7.1. *Let $\{\xi(t), t \in [0,\infty)\}$ be a stochastic process with stationary independent increments and right-continuous paths with $\xi(0) = 0$, $E(\xi(1)) = 1$, $\mathrm{Var}(\xi(1)) = 1$, and*

$$E(e^{\delta|\xi(1)|}) < \infty, \qquad \text{for some } \delta > 0.$$

*Let $\nu$ be a finite measure on $[0,1]$ such that $\nu(\{t\}) = 0$ for all $t \in [0,1]$. Define*

$$\xi_{\theta,\nu}(A) = \xi(\theta\nu(A)), \qquad A \in \mathcal{B}[0,1],$$

*and*

$$\xi_{\theta,\nu}(t) = \xi_{\theta,\nu}([0,t]), \qquad \nu(t) = \nu([0,t]).$$

*Then the family* $\{P(a(\theta)(\xi_{\theta,\nu}(t) - \theta\nu(t))/\theta \in \cdot) : \theta > 0\}$ *satisfies a LDP in* $(D[0,1], \|\cdot\|)$ *with speed* $\frac{a^2(\theta)}{\theta}$ *and rate function*

$$I(\varphi) = \begin{cases} \frac{1}{2} \int_0^1 \left|\frac{d\varphi}{d\nu}(t)\right|^2 \nu(dt), & \text{if } \varphi \ll \nu, \\ +\infty, & \text{otherwise,} \end{cases}$$

*where* $\|\varphi\| := \sup_{t \in [0,1]} |\varphi(t)|$ *for* $\varphi \in D[0,1]$.

PROOF. It suffices to verify the following three conclusions (cf. [2, 31]):

(i) For any $0 < t_1 < t_2 < \cdots < t_k \leq 1$,

$$\left\{ P\left(\frac{a(\theta)}{\theta}(\xi_{\theta,\nu}(t_1) - \theta\nu(t_1), \ldots, \xi_{\theta,\nu}(t_k) - \theta\nu(t_k)) \in \cdot\right), \theta > 0 \right\}$$

satisfies a LDP with the speed $\frac{a^2(\theta)}{\theta}$ and the rate function

$$I_{t_1,\ldots,t_k}(z) = \frac{1}{2} \sum_{i=1}^k \frac{|z_i - z_{i-1}|^2}{\nu((t_{i-1}, t_i])}.$$

(ii) For any $\delta > 0$,

$$(7.1) \quad \lim_{\varepsilon \to 0} \sup_{s \in [0,1]} \limsup_{\theta \to \infty} \frac{a^2(\theta)}{\theta} \log P\left(\sup_{s \leq t \leq s+\varepsilon} |\xi(\theta\nu((s,t])) - \theta\nu((s,t])| \geq \frac{\theta}{a(\theta)} \delta\right) = -\infty.$$

(iii) $I(\varphi) = \sup_{t_1,\ldots,t_k \subset (0,1]} I_{t_1,\ldots,t_k}(\varphi(t_1), \ldots, \varphi(t_k))$.

Since $\xi(t)$ is a random process with stationary and independent increments and the mapping:

$$(z_1, z_2 - z_1, \ldots, z_k - z_{k-1}) \to (z_1, z_2, \ldots, z_k)$$

is continuous in $\mathbb{R}^k$, it is easy to get (i) from the Gärtner–Ellis theorem, the product principle and the contraction principle; (iii) is a consequence of the Cauchy–Schwarz inequality and the submartingale convergence theorem.

Finally, we verify (ii). By Corollary 4 in [24], it is easy to see that there is a universal constant $c > 1$ such that

$$P\left(\sup_{s \leq t \leq s+\varepsilon} |\xi(\theta\nu((s,t])) - \theta\nu((s,t])| \geq \frac{\theta\delta}{a(\theta)}\right)$$

$$\leq cP\left(|\xi(\theta\nu((s,s+\varepsilon])) - \theta\nu((s,s+\varepsilon])| \geq \frac{\theta\delta}{a(\theta)c}\right)$$

$$\leq ce^{-\theta\delta\alpha/(a^2(\theta)c)} E\left(\exp\left\{\frac{1}{a(\theta)}\alpha|\xi(\theta\nu((s,s+\varepsilon])) - \theta\nu((s,s+\varepsilon])|\right\}\right)$$

$$\leq ce^{-\theta\delta\alpha/(a^2(\theta)c)}\left(E\left(\exp\left\{\frac{1}{a(\theta)}\alpha(\xi(1)-1)\right\}\right)^{\theta\nu((s,s+\varepsilon])}\right.$$
$$\left.+E\left(\exp\left\{\frac{-1}{a(\theta)}\alpha(\xi(1)-1)\right\}\right)^{\theta\nu((s,s+\varepsilon])}\right)$$

where $\alpha > 0$ is arbitrary. By the hypotheses, expanding the cumulant yields

$$E\left(\exp\left\{\frac{\pm 1}{a(\theta)}\alpha(\xi(1)-1)\right\}\right) = \exp\left\{\frac{\alpha^2}{2a^2(\theta)} + o(1)\right\}.$$

Therefore

$$\limsup_{\theta\to\infty} \frac{a^2(\theta)}{\theta} \log P\left(\sup_{s\leq t\leq s+\varepsilon} |\xi(\theta\nu((s,t])) - \theta\nu((s,t])| \geq \frac{\theta}{a(\theta)}\delta\right)$$
$$\leq -\sup_{\alpha>0}\left\{\frac{\alpha\delta}{c} - \frac{\alpha^2\nu((s,s+\varepsilon])}{2}\right\} = -\frac{\delta^2}{2c^2\nu((s,s+\varepsilon])}$$

which implies (7.1). □

We now establish the MDP for the Dirichlet process.

THEOREM 7.2. *The family $\{P(a(\theta)(Z_{\theta,\nu}(t) - \nu(t)) \in \cdot) : \theta > 0\}$ satisfies a LDP in $D[0,1]$ with speed $\frac{a^2(\theta)}{\theta}$ and rate function*

$$I_D(\varphi) = \begin{cases} \frac{1}{2}\int_0^1 \left|\frac{d\varphi}{d\nu}(t)\right|^2 \nu(dt), & \text{if } \varphi \ll \nu,\ \varphi(1) = 0, \\ +\infty, & \text{otherwise.} \end{cases}$$

PROOF. Choose $\xi(t)$ in Theorem 7.1 to be the Gamma process $X(t)$. Set

$$Y_{\theta,\nu}(t) = \frac{a(\theta)(X_{\theta,\nu}(t) - \theta\nu(t))}{\theta} - \frac{a(\theta)\nu(t)(X_{\theta,\nu}(1) - \theta)}{\theta}$$
$$= (1-\nu(t))\frac{a(\theta)(X_{\theta,\nu}(t) - \theta\nu(t))}{\theta}$$
$$- \nu(t)\frac{a(\theta)(X_{\theta,\nu}(1) - X_{\theta,\nu}(t) - \theta(1-\nu(t)))}{\theta}.$$

By Theorem 7.1 and the contraction principle, the family $\{P(Y_{\theta,\nu}(t) \in \cdot) : \theta > 0\}$ satisfies a LDP in $D[0,1]$ with speed $\frac{a^2(\theta)}{\theta}$ and rate function

$$\inf\{I(\psi); \psi(t) - \psi(1)\nu(t) = \varphi(t), t \in [0,1]\}$$
$$= \begin{cases} \frac{1}{2}\inf_{\alpha\in\mathbb{R}}\int_0^1 \left|\frac{d\varphi}{d\nu}(t) - \alpha\right|^2 \nu(dt), & \text{if } \varphi \ll \nu,\ \varphi(1) = 0, \\ +\infty, & \text{otherwise,} \end{cases}$$

$$= \begin{cases} \dfrac{1}{2} \int_0^1 \left|\dfrac{d\varphi}{d\nu}(t)\right|^2 \nu(dt), & \text{if } \varphi \ll \nu,\ \varphi(1) = 0, \\ +\infty, & \text{otherwise.} \end{cases}$$

Since

$$|a(\theta)(Z_{\theta,\nu}(t) - \nu(t)) - Y_{\theta,\nu}(t)|$$
$$\leq \left|\dfrac{\theta}{X_{\theta,\nu}(1)} - 1\right| \left(\left|\dfrac{a(\theta)(X_{\theta,\nu}(t) - \theta\nu(t))}{\theta}\right| + \left|\dfrac{a(\theta)\nu(t)(X_{\theta,\nu}(1) - \theta)}{\theta}\right|\right),$$

it follows that for any $\delta > 0$, and $\varepsilon > 0$,

$$P\left(\sup_{t\in[0,1]} |a(\theta)(Z_{\theta,\nu}(t) - \nu(t)) - Y_{\theta,\nu}(t)| > \delta\right)$$
$$\leq P\left(\left|\dfrac{\theta}{X_{\theta,\nu}(1)} - 1\right| > \varepsilon\right) + P\left(\sup_{t\in[0,1]} \left|\dfrac{a(\theta)(X_{\theta,\nu}(t) - \theta\nu(t))}{\theta}\right| > \dfrac{\delta}{2\varepsilon}\right)$$
$$+ P\left(\left|\dfrac{a(\theta)(X_{\theta,\nu}(1) - \theta)}{\theta}\right| > \dfrac{\delta}{2\varepsilon}\right).$$

Now from the LDP of $X_{\theta,\nu}(1)$, one obtains

$$\limsup_{\theta\to\infty} \dfrac{1}{\theta} \log P\left(\left|\dfrac{\theta}{X_{\theta,\nu}(1)} - 1\right| > \varepsilon\right) \leq - \inf_{|1/x-1|>\varepsilon} I(x) < 0,$$

which implies

$$\limsup_{\theta\to\infty} \dfrac{a^2(\theta)}{\theta} \log P\left(\left|\dfrac{\theta}{X_{\theta,\nu}(1)} - 1\right| > \varepsilon\right) = -\infty.$$

From the MDP of $X_{\theta,\nu}$, we have

$$\lim_{\varepsilon\to 0} \limsup_{\theta\to\infty} \dfrac{a^2(\theta)}{\theta} \log P\left(\sup_{t\in[0,1]} \left|\dfrac{a(\theta)(X_{\theta,\nu}(t) - \theta\nu(t))}{\theta}\right| > \dfrac{\delta}{2\varepsilon}\right) = -\infty$$

and

$$\lim_{\varepsilon\to 0} \limsup_{\theta\to\infty} \dfrac{a^2(\theta)}{\theta} \log P\left(\left|\dfrac{a(\theta)(X_{\theta,\nu}(1) - \theta)}{\theta}\right| > \dfrac{\delta}{2\varepsilon}\right) = -\infty.$$

Therefore, for any $\delta > 0$,

$$\limsup_{\theta\to\infty} \dfrac{a^2(\theta)}{\theta} \log P\left(\sup_{t\in[0,1]} |a(\theta)(Z_{\theta,\nu}(t) - \nu(t)) - Y_{\theta,\nu}(t)| > \delta\right) = -\infty;$$

that is, $a(\theta)(Z_{\theta,\nu}(t) - \nu(t))$ is exponentially equivalent to $Y_{\theta,\nu}(t)$. □

**Acknowledgments.** We wish to thank the referees and an Associate Editor for their insightful comments and suggestions. The valuable suggestions of Ian Iscoe are gratefully acknowledged.

DEPARTMENT OF MATHEMATICS AND STATISTICS
MCMASTER UNIVERSITY
HAMILTON, ONTARIO
CANADA L8S 4K1
E-MAIL: shuifeng@mcmaster.ca

SCHOOL OF MATHEMATICS AND STATISTICS
WUHAN UNIVERSITY
WUHAN 430072
P. R. CHINA
E-MAIL: fqgao@whu.edu.cn